\documentclass[10pt]{article}
\usepackage{amsmath, amssymb, amsthm, pb-diagram, euscript, mathrsfs}

\usepackage[T1]{fontenc}
\usepackage[sc]{mathpazo}
\usepackage[usenames]{color}
\usepackage{colortbl}

\DeclareFontFamily{T1}{pzc}{}
\DeclareFontShape{T1}{pzc}{m}{it}{1.8 <-> pzcmi8t}{}
\DeclareMathAlphabet{\mathpzc}{T1}{pzc}{m}{it}

\title{Completely Bounded Characterization of Operator Algebras with Involution}
\author{Nikolay P. Ivankov\\
Max-Planck-Institut f\"ur Mathematik, Bonn}

\theoremstyle{plain}
\newtheorem{prop}{Proposition}[section]

\newtheorem{thm}[prop]{Theorem}

\theoremstyle{definition}
\newtheorem{defn}[prop]{Definition}
\newtheorem{rem}[prop]{Remark}

\theoremstyle{remark}

\newcommand{\A}{\mathcal{A}}

\newcommand{\E}{\mathcal{E}}

\newcommand{\Hc}{\mathcal{H}}
\newcommand{\K}{\mathcal{K}}

\newcommand{\Bb}{\mathbb{B}}
\newcommand{\Cb}{\mathbb{C}}
\newcommand{\Kb}{\mathbb{K}}
\newcommand{\Nb}{\mathbb{N}}

\newcommand{\Gm}{\Gamma}

\newcommand{\Om}{\Omega}

\newcommand{\al}{\alpha}
\newcommand{\be}{\beta}
\newcommand{\gm}{\gamma}

\newcommand{\te}{\theta}

\newcommand{\ro}{\rho}
\newcommand{\sg}{\sigma}

\newcommand{\Id}{\mathrm{Id}}

\newcommand{\diag}{\mathrm{diag}}

\newcommand{\KK}{\mathrm{KK}}
\newcommand{\Ai}[1]{\A^{({#1})}}

\newcommand{\CB}{\mathrm{CB}}

\newcommand{\ox}{\otimes}

\newcommand{\iy}{\infty}

\newcommand{\nix}{\phantom{.}}

\begin{document}

\maketitle
\begin{abstract}
In this paper we study the completely bounded anti-isomorphisms on operator algebras, that work similarly to the involutions with the exception for the property of being completely isometric. We elaborate the Blecher's characterization theorem for operator algebras to make it applicable to the so-called operator $K$-algebras with completely bounded reflexive anti-isomorphism. We also establish a connection of this result with the notion of smooth $C^*$-modules, that play an important role in Mesland's approach to Baaj-Julg picture of $\KK$-theory. 

\end{abstract}

\section{Introduction}

This article is supposed to support the theory developed in \cite{Iva1}. In \cite{Iva1}, in turn, the author tries to enlarge a theoretical basis for the generalization of Kasparov product to the Baaj-Julg picture of $\KK$-theory, which is being developed by Bram Mesland in \cite{Mesl}. 

The main idea of \cite{Mesl} is that under certain conditions called as \emph{transversality} the Kasparov product in $\KK$-theory may be replaced by a simple formula involving the so-called unbounded $\KK$-cycles inroduced by Baaj and Julg. However, in \cite{Mesl}, if $A,B,C$ are $C^*$-algebras and $(E_1, T)$ and $(E_2,S)$ are $(A,B)$ and $(B,C)$ unbounded $\KK$-cycles respectively, then the transversality condition is given in terms of the two concrete unbounded operators $S$ and $T$. In \cite{Iva1} and the author's upcoming thesis we have proposed the way to justify the transversality of the operator $T$ with respect to some class of unbounded $(B,C)$-$\KK$-cycles. For that we have introduced the notion of abstract systems of smooth subalgebas of the $C^*$-algebra $B$ (which differs from the approach of \cite{Mesl} where smooth systems are constructed by means of the unbounded operator $S$), and the $C^k$-algebras in these smooth systems are operator algebras that are either supposed to be given or are constructed in some way. We consider some simplest examples in the end of the article, a more explicit information may be found in \cite{Iva1} and \cite{Mesl}. 

In the case when the algebras are constructed we often encounter the point where one has an operator space which is an algebra, but not operator algebra. Moreover, the construction of smooth algebras in \cite{Mesl} often uses the fact that the involution on the $C^*$-algebra $A$ induces a completely isometric anti-isomorphism on its $C^k$ subalgebras, and, again, this property may be not automatically fulfilled by the construction. Therefore we found a need for a result that would characterize the objects which can be completely boundedly isomorphic to operator algebras with a completely isometric involution. 

In this article we develop such a characterization. It is based on now classical Blecher's characterisation theorem for operator algebras \cite{Ble_cb}, and incorporates an additional involution structure. The result we present here is purely operator algebraic, and may prove itself to be useful in other fields concerning the operator algebra theory.

\section{Preliminaries}

We recall the basic definitions from the theory of operator spaces.
\begin{defn}\label{oper_spc_defn_concr_spce}
A \emph{(concrete) operator space} is a linear subspace of $\Bb(H)$ for some Hilbert space $H$. A \emph{(concrete) operator algebra} is a subalgebra of $\Bb(H)$. The map $f\colon X\to Y$ between two operator spaces is called \emph{completely bounded} ($cb$-map) if there exists a positive constant $C$ such that the natural extensons of $f$ 
\begin{eqnarray*}
f_n\colon M_n(X)& \to & M_n(Y) \\
(x_{ij}) &\mapsto & (f(x_{ij}))
\end{eqnarray*}
has the norm less then $C$ for all $n\in\Nb$. The number $$\|f\|_{cb}:=\sup_n\|f_n\|$$ 
is called the \emph{$cb$-norm of $f$}. The $cb$-map $f$ is called a \emph{$cb$-isomorphism} if it is an isomorphism and its inverse $f^{-1}$ is also completely bounded. If $f$ is a $cb$-isomorphism, such that $\|f\|_{cb}=\|f^{-1}\|_{cb}=1$, then $f$ is called \emph{complete isometry}.
\end{defn}

In what follows we require the operator spaces to be complete with respect to the operator norms on them.

\begin{rem}
One may, of course, have an isomorphism $f\colon X \to Y$ of operator spaces, which is completely bounded, but its inverse $f^{-1}$ is not. Therefore sometimes the term \emph{completely bicontinuous} map is used instead of $cb$-isomorphism. Here, however, we are going to use the term $cb$-isomorphism only in the sense of the Definition \ref{oper_spc_defn_concr_spce}, therefore avoiding the ambiguity.
\end{rem}

\begin{defn}
An operator space $A$ which is also an algebra, such that the multiplication map $m \colon A\times A \to A$ is a completely bounded bilinear map with  $\|m\|_{cb}\le K$ will be called an \emph{operator $K$-algebra} (cf. \cite{BleLeM}). We will use the term \emph{operator pseudoalgebra} when the number $K$ is not specified. The $cb$-homomorphism, $cb$-iomorphism and completely isometric isomorphism of between two operator pseudoalgebras are then algebra homomorphism (isomorphism) which is completely bounded (isometric) as a map between operator spaces.
\end{defn}

\begin{rem}
The notion of operator pseudoalgebra employed in this paper differs from the one given in \cite{Ruan}. In the notation of \cite{BleLeM} these algebras would rather be called operator $1$-algebras.
\end{rem}

There are two important results that give a characterization of operator spaces and operator $K$-algebras respectively. The first one is due to Effros and Ruan.

\begin{thm}[\cite{EfRu}]\label{oper_alg_ruan_oper_space_thm}
Let $X$ be a linear space with a set of matrix norms $\nix_n\|\cdot\|$ on $M_n(X)$, satisfying the properties
\begin{itemize}
\item 
$\nix_{n+m}\|x\oplus y\| = \max\{\nix_n\|x\|,\nix_m\|y\|\}$
\item 
$\nix_n\|\al x \be\| \le \|\al\| \nix_n\|x\| \|\be\|$ 
\end{itemize}
for all $x\in M_n(X)$, $y\in M_m(X)$ and $\al,\be\in M_n(\Cb)$. Then $X$ is completely isometrically isomorphic to a concrete operator space.
\end{thm} 
Thus, we have a characterization of operator spaces up to a complete isometry.

Another theorem is due to Blecher, and it would be the main point of our attention throughout the paper.

\begin{thm}[\cite{Ble_cb},\cite{BleLeM}]\label{oper_alg_thm_cb_charactrization}
Let $A$ be an operator $K$-algebra. Then there exists a (concrete) operator algebra $A'$, which is $cb$-isomorphic to $A$. Moreover, it may be chosen in such a way that if $f\colon A \to A'$ is a $cb$-isomorphism, then $\max\{\|f\|_{cb},\|f^{-1}\|_{cb}\} \le \max \{K^{-1},2K\}$.
\end{thm}

\begin{rem}
Obviously, all the concrete operator algebras are operator $1$-algebras. The converse in general is not true. To have an operator $1$-algebra being completely isometrically isomorphic to a concrete operator algebra, one has to add an assumption that $A$ possesses a contractive approximate unit (cf. \cite{Ruan})
\end{rem}

\section{Involution}

Recall that an involution on a Banach algebra $A$ is an isometric anti-isomorphism $\ast \colon A \to A$, $\ast\colon a\mapsto a^*$ such that $a^{\ast\ast}=a$. 

Thus, if we want to specialize this notion for the case of operator algebras, we should first give a definition of a $cb$-anti-isomorphism.

\begin{defn}
Let $A$ be an operator pseudoalgebra. Then an anti-homomorphism $f\colon A \to B$ will be called \emph{$cb$-anti-homomorphism} if there exists a positive number $C$ such that
$$
\nix_n\|(f(a_{ji}))_{ij}| \le C\nix_n\|(a_{ij})_{ij}\|
$$
for all $(a_{ij})_{ij}\in M_n(A)$. If $f$ is anti-isomorphic, and its inverse $f^{-1}$ is also a $cb$-anti-homomorphism, then $f$ will be called a \emph{$cb$-anti-isomorphism}. One may analogously define completely isometric anti-isomorphisms. 
\end{defn}

\begin{rem}
Observe that, unlike the case of homomorphisms, we have to add a transposition in matrix algebras to the definition of $cb$-anti-homomorphisms. This makes the notion of $cb$-anti-homomorphism much more subtle then the one of $cb$-homomorphism. It seems, although the author doesn't have a concrete example for now, that even for a general (concrete) operator algebra $A$ there would not be any $cb$-anti-isomorphisms of $A$ onto itself. However, as we have indicated in the introduction and will also see in the next section, the algebras having $cb$-anti-isomorphisms may often appear in applications.
\end{rem}

\begin{defn}
A $cb$-anti-isomorphism $f\colon A\to A$ such that $f^2 = \Id_A$ would be called an \emph{(operator algebra) pseudo-involution} on $A$. If, in addition, $f$ is completely isometric then it will be called an \emph{(operator algebra) involution}. An operator algebra possessing an involution will be called \emph{involutive}.
\end{defn}

We are going to show that any pseudo-involution may in some sense be "updated" to become an involution.

\begin{prop}\label{oper_alg_prop_pseudo_involution}
Let $A$ be an operator $K$-algebra with a pseudo-involution $f$. Then there is an operator pseudoalgebra $B$ and a $cb$-isomorphism $\sg\colon A\to B$, such that $\sg f\sg^{-1}$ is an involution on $B$.
\begin{proof}
Let $B=A$ as a algebras We define matrix norms on $B$ as 
$$
\nix_n\|(a_{ij})_{ij}\|_B = \max\{\nix_n\|(a_{ij})_{ij}\|_A,\,\nix_n\|(f(a_{ji})_{ij})\|_A\}
$$

The space $B$ endowed with this norms is an operator pseudoalgebra. Indeed, we have that
\[
\begin{split}
& \nix_{n+m}\|(a_{ij}\oplus b_{kl})\|_B = \\
& = \max \{ \max\{\nix_n\|(a_{ij})_{ij}\|_A,\,\nix_n\|(f(a_{ji}))_{ij}\|_A\}, \, \max\{\nix_m\|(b_{kl})_{lk}\|_A,\,\nix_m\|(f(b_{kl}))_{lk}\|_A\} \}\\
& =
\max \{\max\{\nix_n\|(a_{ij})_{ij}\|_A,\nix_m\|(b_{kl})_{lk}\|_A\},\,
\max \{\nix_n\|(f(a_{ji}))_{ij}\|_A, \nix_m\|(f(b_{kl}))_{lk}\|_A\}\} \\
&= \max \{\nix_{n+m}\|(a_{ij})_{ij}\oplus (b_{kl})_{lk} \|_A, \,
\nix_{n+m}\|(f(a_{ji}))_{ij}\oplus(f(b_{kl}))_{lk}\|_A\} \\
& = 
\max \{\nix_n\|(a_{ij})\|_B, \nix_m\|(b_{kl})\|_B\}
\end{split}
\]
and 
\[
\begin{split}
\nix_n\|\al(a_{ij})\be\|_B & = \max\{\nix_n\|\al(a_{ij})_{ij}\be\|_A, \, \nix_n\|\be^\intercal(f(a_{ji}))_{ij}\al^\intercal\|_A\} \\
&\le
\max\{\|\al\|\nix_n\|(a_{ij})_{ij}\|_A\|\be\|,\, \|\be^\intercal\|\nix_n\|(f(a_{ji}))_{ij}\|_A \|\al^\intercal\|\} \\
& =
\|\al\| \|\be\| \max\{\nix_n\|(a_{ij})_{ij}\|_A, \, \nix_n\|(f(a_{ji}))_{ij}\|_A\}
\end{split}
\]

here we use the fact that $\al$ and $\be$ are scalar matrices. Thus, $B$ is an operator space. To prove that it is a pseudoalgebra, observe that
\[
\begin{split}
\nix_n\|(a_{ij})(b_{kl})\|_B & = \max\{\nix_n\|(a_{ij})(b_{kl})\|_A, \, 
\nix_n \|f_n((a_{ji})(b_{kl}))\|_A\} \\ 
&\le \max \{\nix_n\|(a_{ij})(b_{kl})\|_A, \|f\|_{cb} \nix_n\|(a_{ij})(b_{kl})\|_A\} \\
&\le \|f\|_{cb} K\nix_n\|(a_{ij})\|_A\nix_n\|(b_{kl})\|_A \\
&\le \|f\|_{cb} K \cdot\|f\|_{cb}\max\{\nix_n\|(a_{ij})_{ij}\|_A, \, \nix_n\|(f(a_{ji}))_{ij}\|_A\}
\cdot \\
&\cdot \|f\|_{cb}\max\{\nix_n\|(b_{kl})_{kl}\|_A, \, \nix_n\|(f(b_{lk}))_{kl}\|_A\} \\
& = \|f\|_{cb}^3 K \nix_n\|(a_{ij})\|_B \nix_n\|(b_{kl})\|_B
\end{split}
\]
Here we use the fact that since $f^2=1$ we have that $\|f\|_{cb}\ge 1$.

Since $f$ is a $cb$-anti-isomorphism and $f^2=1$, we have that
$$ 
\|f\|^{-1}_{cb}\nix_n\|\cdot\|_A \le \nix_n\|\cdot\|_B \le \|f\|_{cb}\nix_n\|\cdot\|_A
$$
so the algebras $A$ and $B$ are $cb$-isomorphic. Denote this isomorphism by $\sg$. By the construction $(\sg f \sg^{-1})=\Id_B$. Now, for $(a_{ij})\in M_n(A)$ we have that
\[
\begin{split}
\nix_n\|\sg f \sg^{-1}(\sg(a_{ij}))\|_B &= \nix_n\|\sg f(a_{ij})\|_B \\ 
& = \max\{\nix_n\|f(a_{ij})\|_A,\, \nix_n\|f^2(a_{ij})\|_A\} \\
& = \max\{\nix_n\|f(a_{ij})\|_A,\nix_n\|(a_{ij})\|_A\} = \\
& \nix_n \|\sg(a_{ij})\|_B 
\end{split}
\]
Since $\sg$ is a $cb$-isomorphism, all the elements of $M_n(B)$ have the form $\sg(a_{ij})$. This last observation settles the proof.
\end{proof}
\end{prop}

\begin{rem}
Observe that since $f$ was an anti-isomorphism, we were not able to define $\sg$ as just $\sg\colon a \mapsto a \oplus f(a)$, since in this case $\sg(ab) = ab \oplus f(ba)$.
\end{rem}

The result \ref{oper_alg_prop_pseudo_involution} gives us only an operator pseudoalgebra with involution. However, a closer look to the Theorem \ref{oper_alg_thm_cb_charactrization} lets us extend this result, making $B$ into a (concrete) operator algebra with involution.

In order to do this, we recall the construction from \cite{Ble_cb}. Let $\Gm$ be the set, $n\colon \Gm \to \Nb$, $\gm \mapsto n_\gm$ be a function. Let $\Lambda$ be a set of formal symbols (variables) $x^\gm_{ij}$, one variable for each $\gm\in\Gm$ and  each $1\le i,j \le n_\gm$. Denote by $\Phi$ a free associative algebra on $\Lambda$. In this case $\Phi$ consists of polynomials in the non-commuting variables with no constant term. One then defines a norm on $M_n(\Phi)$ by
\begin{equation}\label{oper_alg_eqn_involution_reprs}
\|(u_{ij})\|_{\Lambda} := \sup_\pi(\|(\pi(u_{ij}))\|)
\end{equation}
where $\pi$ goes through all the representations of $\Phi$ on a separable Hilbert space satisfying the condition $\|(\pi(x^\gm_{ij}))_{ij}\|\le 1$ for all $\gm$, where the latter matrix is indexed on rows by $i$ and on columns by $j$ for all $1\le i,j \le n_\gm$.

It is then shown in \cite{Ble_cb} that the map defined above indeed defines a norm on $M_n(\Phi)$ and that $\Phi$ becomes an operator algebra with respect to these operator norms.

In the proof of the characterization theorem the set $\Gm$ is taken to be the collection of $n\times n$ matrices $\gm=(a_{ij})$ with entries in $A$ such that $\|\gm\|=\frac{1}{2K}$, where $K$ is a $cb$-norm of the multiplication in $A$. Then one takes $\Lambda$ to be the collection of entries of these matrices $x^\gm_{ij}:=a_{ij}$, regarded as formal symbols indexed by $\gm$ and $i,j$, not identifying "equal" entries for different indexes. After that there is defined a 
map 
\begin{eqnarray*}
\te\colon \Phi &\to& A \\
x^\gm_{ij} &\mapsto& \gm_{ij}
\end{eqnarray*} 
which is then extended to general polynomials. It is proved in \cite{Ble_cb} that $\te$ is a completely contractive homomorphism. One then puts $B:= \Phi/\ker(\te)$, which is an operator algebra subject to the quotient operator norm, and is $cb$-isomorphic to $A$.

Now let the pseudoalgebra $A$ be involutive. Observe that since the involution on $A$ is completely isometric, we have that $\nix_n\|(a_{ij})^*\|=\nix_n\|(a_{ij})\|$, and thus $(a_{ij})^*\in\Gm$. Hence we have that $a_{ij}^*\in \Lambda$. This observation makes us able to define an involution the following way. On $\Phi$ we set
$$
(x^{\gm_1}_{i_1j_1} x^{\gm_2}_{i_2j_2} \dots x^{\gm_k}_{i_kj_k})^* := 
(x^{\gm_k}_{i_k j_k})^* (x^{\gm_{k-1}}_{i_{k-1}j_{k-1}})^* \dots (x^{\gm_1}_{i_1j_1})^*
$$ 
on the monomials, and then extend this to the whole $\Phi$. Analogously, on $M_n(\Phi)$ we set 
$$((P_{ij})_{ij})^*=(P_{ji}^*)_{ij}$$

By the construction we have that $\te((P_{ij})^*)=\te((P_{ji}^*)_{ij})$. Consequently, let $\pi\colon \Phi \to \Bb(\Hc)$ be a representation of $\Phi$ satisfying the condition $\|(\pi(x^\gm_{ij}))\|\le 1$ for all $(x^\gm_{ij})_{ij}$. Denote this set by $\Theta$. We may define a representation $\pi' \colon \Phi \to \Bb(\Hc)$ by setting $\pi'((P_{ij})^*):=(\pi(P_{ij}))^*$, where the latter involution is given by the one on the $\Bb(\Hc^{n_\gm})$. By the construction, we have that
$$
\nix_{n_\gm}\|\pi'(x^\gm_{ij})\|=\nix_{n_\gm}\|(\pi((x^\gm_{ij})^*))^*\| =
\nix_{n_\gm}\|\pi((x^\gm_{ij})^*)\| \le 1
$$ 
for all $(x^\gm_{ij})_{ij}$ since $(x^\gm_{ij})^*\in\Gm$, and so $\pi'\in \Theta$. Therefore we have that 
\[
\begin{split}
\|(P_{ij})\|_\Lambda &= \sup_{\pi\in \Theta}(\|(\pi(P_{ij}))\|) \\
& = \sup_{\pi'\in\Theta}\|((\pi'(P_{ij})^*))^*\| \\
& =  \sup_{\pi'\in\Theta}\|\pi'(P_{ij})^*\| \\
& =
\|(P_{ij})^*\|_\Lambda
\end{split}
\]

Hence we obtain that the map $\te$ respects the involution, and thus the anti-isomorphism induced on $B$ by the involution on $\Phi$ preserves the operator norms.

Combining this observations with Proposition \ref{oper_alg_prop_pseudo_involution} we have the following

\begin{thm}\label{oper_alg_thm_making_involution}
Let $A$ be an operator pseudoalgebra and let $f$ be a pseudo-involution on $A$. Then there is a $cb$-isomorphism $\ro \colon A \to B$, such that the map $\ro f \ro^{-1}$ is an (operator algebra) involution on $B$.
\begin{proof}
Put $\ro = \te \sg$.
\end{proof}
\end{thm}
\begin{rem}
We may also estimate the $cb$-norm of $\ro$. Indeed, the map $\sg$ has the $cb$-norm $\|f\|_{cb}$, and gives us a pseudoalgebra $B'$ with the $cb$-norm of multiplication bounded by $\|f\|_{cb}^3K$. Thus, for $K \ge 1$ the estimation from \cite{Ble_cb} shows us that the map $\te$ has a $cb$-norm $\le  2 \|f\|_{cb}^3K$. Hence, $\|\ro\|_{cb} \le 2\|f\|_{cb}^4K $.
\end{rem}

\section{Application $C^1$-Modules} 

In this section we are going to show the relation of the construction of involutive operator algebras to the notion of smooth modules as they are defined in \cite{Mesl} and \cite{Iva1}. We will give here a simplified definition of smooth algebras and modules. For a more descriptive picture, see \cite{Mesl}, \cite{Iva1}.

Let $A$ be a $C^*$-algebra $E$ be a Hilbert $C^*$-module over $A$, and let $\A$ be an operator algebra, which is isomorphic to a pre-$C^*$-subalgebra of $A$ abusively denoted by $\A$. We define a $C^1$ structure on $E$ with respect to $\A$ by choosing a countable approximate unit $u_n = \sum_{j=1}^{n} x_j \ox x_j$ on $\Kb_A(E)$ with a property that 
$$
\|(\langle x_j, x_k\rangle)_{jk}\|_{1,D} \le C
$$ 
Then, a pre-$C^1$-module over the $C^1$ algebra $\A$ is defined as 
$$
\E = \{e \in E \mid \langle x_j,e \rangle \in \A,\, \|\sum_{j=1}^\iy\langle e, x_j \rangle\|_1 < \iy\}
$$
and it is a $C^1$-module if it satisfies the Kasparov stabilization property.

Now if the involution on $A$ induces an operator algebra involution on $\A$, then the space $\Kb_{\A}(\E,\A)$ is completely isometrically isomorphic to  $\E$. Thus, there is a well-defined inner product on $\E$, which is a restriction of the inner product on $\E$. The existence of this product then allows use to construct canonically the algebra $\CB^*_\A\E$ of completely bounded $\A$-linear involutive of operators on $\E$.

Suppose now that there is another operator algebra $\A'$ with the same properties as $\A$, such that $\A \hookrightarrow \A'$ as pre-$C^*$-algebras, and the inclusion map induces a completely bounded injective homomorphism of corresponding operator algebras. Then, by the construction, the smooth structure on $E$ with respect to $\A$ will automatically be a smooth structure with respect to $\A'$, and we obtain a completely bounded inclusion $\E \to \E'$, where $\E'$ is obtained form the approximate unit $u_k$ analogously to $\E$. This observation also allows us to transfer additional structures which are involved in the construction of $\KK$-product from $\E$ to $\E'$.

This construction may then become an intermediate step in the construction of $KK$-product in the Baaj-Julg picture of $\KK$-theory. We briefly describe the simplest case.

Let $A$, $B$ be $C^*$-algebras and $(E,D)$ be an unbounded $(A,B)$-$\KK$-cycle on  (see, \cite{BaJu}, \cite{Bla} for definition), and we suppose for a moment that $D$ is selfadjoint. Denote
$$
\Ai{1}_D := \{a\in A \mid [D,a] \text{ extends to an element of } \CB^*_A(E)\}
$$ 
Here we use the graded commutator. By the definition of an unbounded $\KK$-cycle the algebra $\Ai{1}_D$ is dense in $A$.

We introduce a representation of $\Ai{1}_D$ by setting
$$
\pi^1_D(a) = \begin{pmatrix}
a & 0 \\
[D,a] & a
\end{pmatrix}
$$
with the operator norm $\|\cdot\|_{1,D}$ defined by this representation.
So, by definition $\Ai{1}_D$ is a concrete operator algebra. The involution on $A$ induces an operator algebra pseudoinvolution on $\Ai{1}_D$. Indeed,
\[
\begin{split}
\|a^*\|_{1,D} &= \left\|\begin{pmatrix} a^* & 0 \\ [D;a^*] & a^*\end{pmatrix}\right\| \\
& =  \left\|\begin{pmatrix} a & \pm[D;a] \\ 0 & a\end{pmatrix}\right\| \\
&= \left\|
\begin{pmatrix} 0 & \Id_E \\ \Id_E & 0\end{pmatrix}
\begin{pmatrix} a & \pm[D;a] \\ 0 & a\end{pmatrix}
\begin{pmatrix} 0 & \Id_E \\ \Id_E & 0\end{pmatrix}
\right\| \\
&=\left\|\begin{pmatrix} a & 0 \\ \pm[D;a] & a\end{pmatrix}\right\|
\end{split}
\]
where the sign is $+$ when there the cycle $(E,D)$ is even and $-$ when it is odd. We will stick for now to the even case. Then the involution will obviously be isometric. To show that the involution is completely isometric, observe that for each $m$ there is a "permutation" unitary $U_m$, such that
\[
U_m (\pi^1_D(a_{jk}))_{jk} U_m^{-1} = 
\begin{pmatrix}
(a_{jk}) & 0 \\
[\diag_m(D); (a_{jk})] & (a_{jk})
\end{pmatrix}
= \pi^1_{\diag_m(D)}((a_{jk}))
\]
and so we can use the previous observation.
 
In \cite{Iva1} we construct a kind of "universal" $C^1$-subalgebra for a separable $C^*$-algebra $A$. More precisely for any set of $C^*$-algebras $\Xi$ we construct an operator algebra $\Ai{1}$, such that for any $C^*$-algebra $B_\xi\in\Xi$ there is a set $\Om_\xi$ of unbounded $(A,B)$-$\KK$-cycles $(E_{\xi,\omega},D_{\xi,\omega})$, for which the map $\Om_\xi \to \KK_0(A,B_{\xi})$ given by 
$$(E_{\xi,\omega},D_{\xi,\omega}) \mapsto [(E_{\xi,\omega}, D_{\xi,\omega}(1+D_{\xi,\omega}^*D_{\xi,\omega})^{-\frac{1}{2}})]$$   
and there is a $cb$-inclusion $\Ai{1}\hookrightarrow \Ai{1}_D$ which preserves the involution. The existence of such inclusion guarantees us that if an approximate unit $u_k$ defines a $C^1$ structure on an Hilbert $C^*$ $A$-module $E$ with respect to $\Ai{1}$, then so does it for all $\Ai{1}_{D_{\xi, \omega}}$. 

The idea of the construction is follows. We take a set of representatives $(E_{\xi,\omega}, F_{\xi,\omega})$ of the elements of $\KK_0(A,B_\xi)$, fix a total system $\{a_i\}$ on $A$ and construct unbounded regular selfadjoint operators $(E_{\xi,\omega},D_{\xi,\omega})$ such that 
$$
[(E_{\xi,\omega}, D_{\xi,\omega}(1+D_{\xi,\omega}^*D_{\xi,\omega})^{-\frac{1}{2}})] = [(E_{\xi,\omega}, F_{\xi,\omega})]
$$ 
We also define them in such a way that $\|a_i\|_{D_{\xi,\omega}}\le C_j$ for all $D_{\xi,\omega}$ for all the elements of the chosen total system on $A$. Then we define the algebra 
$$
\Ai{1} = \{a\in A \mid \sup_{\xi,\omega}\|a\|_{1,D_{\xi,\omega}} < \iy\}
$$
and with the collection of matrix norms on it defined as
$$
\nix_m \|(a_{kl})\|_{1} := \sup_{\xi,\omega} \nix_m \|(a_{kl})\|_{1,D_{\xi,\omega}}
$$
Since all the elements the total system $\{a_i\}$ lay in $\Ai{1}$, the algebra $\Ai{1}$ is dense in $A$. It is also shown in \cite{Iva1} that $\Ai{1}_{D_{\xi,\omega}}$ and $\Ai{1}$ are stable under holomorphic functional calculus on $A$ and therefore have the same $\K$-theory as $A$.

It is easy to check then that $\Ai{1}$ is then an operator $1$-algebra with a completely isometric involution. We also have that there is a completely contractive inclusion $\Ai{1}\hookrightarrow\Ai{1}_{D_{\xi,\omega}}$. However, in case when $A$ is nonunital, $\Ai{1}$ may be not isomorphic to a concrete operator algebra. Therefore, in order to make $\Ai{1}$ into an involutive operator algebra, we need to use Theorem \ref{oper_alg_thm_making_involution}.

Observe that if we would like to incorporate the odd modules in the picture, the involution on the algebras $\Ai{1}_{D_{\xi,\omega}}$ will not necessarily be completely isometric any more, and we shall need to use the Theorem \ref{oper_alg_thm_making_involution} to obtain involutive operator algebras.

Another example where the Theorem \ref{oper_alg_thm_making_involution} may become useful arises when one considers \emph{almost selfadjoint} unbounded operators instead of just selfadjoint ones. Let $D$ be a selfadjoint regular operator on a Hilbert $C^*$-$B$-module $E$ and suppose that $b\in \CB^*_B(E)$ and is even, but in general we do not demand the selfasjointness of $b$. We construct an operator algebra $\Ai{1}_{D+b}$ analogously to $\Ai{1}_D$. But now also in the case when we consider even unbounded $\KK$-cycles the involution on $\Ai{1}_{D+b}$, although completely bounded, may be not isometric, since
\[
\begin{split}
\|a^*\|_{1,D+b} &= 
\left\|
\begin{pmatrix}
a^* & 0 \\
[D+b,a^*] & a^*
\end{pmatrix}
\right\| \\
& =
\left\|
\begin{pmatrix}
a & [D+b^*,a] \\
0 & a
\end{pmatrix} 
\right\| \\
& =
\left\|
\begin{pmatrix}
a & 0 \\
[D+b^*,a] & a
\end{pmatrix} 
\right\| 
\end{split}
\]
and the latter norm should not in general be equal to $\|a\|_{1,D}$. The Proposition \ref{oper_alg_prop_pseudo_involution} then tells us that there is a (canonical) way to associate an involutive operator algebra to $\Ai{1}_{D+b}$ and therefore simplify the consequent calculations. 

In the in the theory developed in \cite{Mesl} one may encounter other examples of operator algebras in which the involution should not necessarily be completely isometric, but only completely bounded. We have already mentioned one of the most important of them: these are the algebras of the form $CB^*_{\A}(\E)$ of completely bounded involutive $\A$ operators over $C^1$-module $\E$ and their involutive subalgebras. The latter ones with operator norms induced by the norm on $CB^*_{\A}(\E)$ are used for the definition of subsequent $C^1$-modules, which, in turn, are used for further construction of unbounded Kasparov product.

Finally, it should be noted that in \cite{Iva1} and \cite{Mesl} there are considered higher orders of smoothness of algebras and modules, and the results we have presented in this paper may be also applicable to these cases.


\begin{thebibliography}{100}

\bibitem{BaJu} Saad Baaj and Pierre Julg, \emph{Th\'eorie Bivariante de
Kasparov et Op\'erateurs non Bornes dans les $C^*$-Modules Hilbertiens}. C.R. Acad Sci. Paris, No. 296 (1983), Ser. I, pp. 875-878




\bibitem{Bla} Bruice Blackadar. \emph{K-Theory for Operator Algebras.} Springer-Verlag New York Inc., 1986



\bibitem{Ble_cb} David P. Blecher, \emph{A completely Bounded Characterisation of Operator ALgebras}, Mathematische Annalen, No. 303 (1995), pp. 227-239.



\bibitem{BleLeM} David P. Blecher, Christian Le Merdy, \emph{Operator Algebras and Their Modules - An Operator Space Approach}, Oxford Univ. Press, 2004.










\bibitem{EfRu} E. Effros, Zhong-Jin Ruan, \emph{On the Abstract Characterization of Operator Spaces}, Proc. Amer. Math. Soc., No. 119 (1993), pp. 579-584

\bibitem{Iva1} N. P. Ivankov, \emph{Noncommutative Fr{\'e}chet Spaces and Unbounded Bivariant K-Theory}, arXiv:1103.4528v1 [math.KT]





\bibitem{Mesl} Bram Mesland, \emph{Bivariant K-Theory of Groupoids and the 
Noncommutative Geometry of Limit Sets}, arXiv:0904.4383v2 [math.KT]








\bibitem{Ruan} Zhong-Jin Ruan, \emph{A characterization of nonunital operator 
algebras}, Proc. of AMS, Vol. 121, No. 1 (May, 1994), pp. 193-198  





\end{thebibliography}
\end{document}